\documentclass[12pt,a4paper]{amsart}
\usepackage[utf8]{inputenc}	
\usepackage[T1]{fontenc}
\usepackage[in,headings]{fullpage}
\usepackage{amsbsy, amscd, amsfonts, amsmath, amsrefs, amssymb, amsthm}
\usepackage{graphicx}
\usepackage{float}
\usepackage{color}   
\usepackage[colorlinks=true,linkcolor=blue,citecolor=blue,urlcolor=blue]{hyperref}

\newtheorem{theorem}{Theorem}

\newtheorem{corollary}[theorem]{Corollary}
\newtheorem{proposition}[theorem]{Proposition}

\theoremstyle{definition}

\newtheorem{remark}[theorem]{Remark}

\theoremstyle{remark}

\newcommand{\C}{\mathbf{C}}

\newcommand{\R}{\mathbf{R}}
\newcommand{\N}{\mathbf{N}}

\renewcommand{\Re}{\mathop{\mathrm{Re}}\nolimits}

\newcommand{\Rzeta}{\mathop{\mathcal R }\nolimits}

\newfont{\cmbsy}{cmbsy10}
\newfont{\cmmib}{cmmib10}

\DeclareMathOperator*{\Res}{Res}

\begin{document}

\title[Entire function defined by Riemann]
{An entire function defined by Riemann}
\author[Arias de Reyna]{J. Arias de Reyna}
\address{%
Universidad de Sevilla \\ 
Facultad de Matem\'aticas \\ 
c/Tarfia, sn \\ 
41012-Sevilla \\ 
Spain.} 




\email{arias@us.es, ariasdereyna1947@gmail.com}


\begin{abstract}
In one of the sheets in Riemann's Nachlass he defines an entire function
and connect it with his  zeta function. As in many pages in his Nachlass, 
Riemann is not giving complete proofs.

He obtains an $L^\infty$ function   whose Fourier transform vanish at the 
real values $\gamma$ with $\zeta(\frac12+i\gamma)=0$.  We give proofs of 
Riemann formulas.

This is an integral representation of the zeta function different from the known ones. I believe this is the first time it has been published. 
\end{abstract}

\maketitle

\section{Introduction}

After the publication of Riemann's complete works, it was suspected that his posthumous writings contained an expression of $\Xi(t)$. Thanks to the contributions of the librarian Distel and of Bessel Hagen, the Riemann-Siegel formula and the integral expression $\Rzeta(s)$ were discovered. They were studied by Siegel in \cite{Siegel}. 
Here, we consider a different integral representation that we have found in Riemann's posthumous papers.\footnote{I have been able to see Riemann's posthumous writings on number theory thanks to W. Gabcke who sent them to me with permission from the library in Gotingen.   The integral expression under consideration is contained in page 00000093.jpg} 

As usual in Riemann we only find a set of equations. I  give a proof of his formulas and some notice of my interpretation of these formulas. 
Riemann defines a bounded function on $\R$ whose Fourier transform vanishes at the real values $\gamma$ with $\zeta(\frac12+i\gamma)=0$. Riemann defines an entire function
\begin{equation}
\chi(z):=\int_L\frac{e^{izw}}{e^{2\pi w^2}-1}\,dw,
\end{equation}
where $L$ is a line parallel to the real axis at a height $0<b<2^{-1/2}$. It has the power series expansion
\begin{equation}
\chi(z)=\frac{z}{2}+\frac{1}{\sqrt{2}}\sum_{n=0}^\infty \frac{\zeta(n+\frac12)}{n!}\Bigl(-\frac{z^2}{8\pi}
\Bigr)^n,
\end{equation}
and satisfies the functional equation $\chi(-z)=\chi(z)-z$. And is connected with $\zeta(s)$ by 
\begin{equation}
\int_0^\infty \chi(x)x^{2s-1}\,dx=\Gamma(2s)\zeta(s+\tfrac12)\cos\tfrac{\pi(2s+1)}{4},
\qquad \Re(s)>0.
\end{equation}

\section{Definition of Riemann's chi function.}

\begin{proposition}
The integral
\begin{equation}\label{e1}
\chi(z):=\int_L\frac{e^{izw}}{e^{2\pi w^2}-1}\,dw,
\end{equation}
where $L$ is a line parallel to the real axis at a height $0<b<2^{-1/2}$
define an entire function.
\end{proposition}

\begin{proof}
For a fixed $z$ the integrand is meromorphic in  $w$ with simple poles at 
$\pm \sqrt{n}\;e^{\pm \frac{\pi i}{4}}$ with $n\in \N$ and a double pole at $w=0$. 

Let $R>0$ and $0<\varepsilon<2^{-1/2}$ be positive real numbers. 
For $|z|<R$ and $w=x+ib$ with $\varepsilon<b<2^{-1/2}-\varepsilon$ and $x$ real,  we have
\[\Bigl|\frac{e^{izw}}{e^{2\pi w^2}-1}\Bigr|\le C e^{-\pi x^2}\]
for an adequate constant $C=C(R,\varepsilon)$. Therefore $\chi(z)$ is well defined and holomorphic in $|z|<R$.
\end{proof}

\section{Calculating the integral by residues.}

Riemann first obtains a series for his function $\chi(z)$.
\begin{proposition}
For $|\arg(z)|<\frac{\pi}{4}$ we have
\begin{equation}\label{e5}
\chi(z)=\frac12\sum_{n=1}^\infty \Bigl(\frac{e^{-z\sqrt{n}e^{-\pi i/4}}}{\sqrt{n}}
e^{\pi i/4}
+\frac{e^{-z\sqrt{n}e^{\pi i/4}}}{\sqrt{n}}e^{-\pi i/4}\Bigr).
\end{equation}
In particular, for $z=x>0$ real 
\begin{equation}\label{e6}
\chi(x)=\sum_{n=1}^\infty\frac{e^{-x\sqrt{n/2}}}{\sqrt{n}}\sin\Bigl(\,\frac{\pi}{4}-
x\,\sqrt{\frac{n}{2}}\;\Bigr).
\end{equation}
\end{proposition}

\begin{proof}
Let $N$ be a non-negative integer and define $b$ by  
$b\sqrt{2}=N+\frac12$. Let $L_N$ be the 
path consisting of the line $x+ib$.  This line is chosen so that it passes 
between successive poles of the integrand defining $\chi(z)$. 

By Cauchy's Theorem 
\begin{multline*}
\int_{L_0}\frac{e^{izw}}{e^{2\pi w^2}-1}\,dw- 
\int_{L_N}\frac{e^{izw}}{e^{2\pi w^2}-1}\,dw\\=
2\pi i \sum_{n=1}^N\Bigl(\Res_{w=-\sqrt{n}e^{-\pi i/4}}\frac{e^{izw}}{e^{2\pi w^2}-1}
+\Res_{w=\sqrt{n}e^{\pi i/4}}\frac{e^{izw}}{e^{2\pi w^2}-1}\Bigr).
\end{multline*}
After some computations
\[\chi(z)=\frac12\Bigl(\sum_{n=1}^N\frac{e^{-z\sqrt{n}e^{\pi i/4}}}{\sqrt{n}}e^{-\pi i/4}+ \frac{e^{-z\sqrt{n}e^{-\pi i/4}}}{\sqrt{n}}
e^{\pi i/4}\Bigr)+
\int_{L_N}\frac{e^{izw}}{e^{2\pi w^2}-1}\,dw.\]

For a fixed $z=a>0$ we show that the integral in $L_N$ converges to $0$ when $N\to\infty$.
In fact, we have
\[U_N:=\int_{L_N}\frac{e^{iaw}}{e^{2\pi w^2}-1}\,dw=
\int_{-\infty}^{+\infty}\frac{e^{-ab}e^{iax}\,dx}{e^{2\pi(x^2-b^2)+4\pi i b x}-1}.\]
Therefore,
\begin{align*}
|U_N|&\le e^{-ab}\int_{-\infty}^{+\infty}\frac{dx}{|e^{2\pi(x^2-b^2)}-e^{-4\pi i b x}|}\\
&=  be^{-ab}\int_{-\infty}^{+\infty}\frac{dx}{|e^{2\pi b^2(x^2-1)}-e^{-4\pi i b^2 x}|}.
\end{align*}
Recall that $2b^2=(N+\frac12)^2$. It can be shown (we will show this a little further on)  that there is a constant $C>0$ and 
an integer $N_0$  such 
that for $N\ge N_0$ and all $x$
\begin{equation}\label{e7}
|e^{2\pi b^2(x^2-1)}-e^{-4\pi i b^2 x}|> C b^{-2}.
\end{equation}
Assuming this, we have for $N\ge N_0$
\begin{align*}
|U_N|&\le be^{-ab}\frac{4}{C}b^2+be^{-ab}\int_{|x|>2}\frac{dx}{|e^{2\pi b^2(x^2-1)}-e^{-4\pi i b^2 x}|}\\
&\ll b^3e^{-ab}+ be^{-ab}\int_2^{+\infty}e^{-2\pi b^2(x^2-1)}\,dx .
\end{align*}
Since this tends to $0$ for $b\to+\infty$, we have $\lim_{N\to\infty}U_N=0$. 

To prove \eqref{e7} observe that \eqref{e7} is  true, with $C=1$ when 
$|e^{2\pi b^2(x^2-1)}-1|\ge b^{-2}$. Therefore,
we may assume that $|e^{2\pi b^2(x^2-1)}-1|< b^{-2}$. In this case, for $n\ge N_0$, we have
\[-\frac{1}{2}b^{-2}<\log(1-b^{-2})<2\pi b^2(x^2-1)<\log(1+b^{-2})<2b^{-2}.\]
This implies that $|x-1|\le cb^{-4}$ or $|x+1|\le cb^{-4}$ for some absolute constant
$c>0$. This implies that $e^{4\pi i b^2x}$ is very near $i$ or $-i$, for example, in the first case $x=1+\theta b^{-4}$, for some real $\theta$ with $|\theta|\le c$ 
\[e^{4\pi i b^2x}=e^{2\pi i(N^2+N+\frac14)(1+\theta b^{-4})}=e^{\frac{\pi i}{2}+4\pi i \theta b^{-2}},\]
so that for a new absolute  constant $c>0$ we have $|e^{4\pi i b^2x}-i|\le cb^{-2}$.
It follows that 
\begin{multline*}
|e^{2\pi b^2(x^2-1)}-e^{-4\pi i b^2 x}|=
|1-i+(e^{2\pi b^2(x^2-1)}-1)-(e^{-4\pi i b^2 x}-i)|\\
\ge \sqrt{2}-b^{-2}-cb^{-2}\ge 
\frac12>C b^{-2}.
\end{multline*}

This proof \eqref{e5} when $z$ is real and positive. But it is easy to see that
both members are well defined and holomorphic in $|\arg(z)|<\pi/4$. So the 
equality extends to this range.  Equation \eqref{e6} is obtained from \eqref{e5}
by using the definition of the $\sin$ function. It is also valid in the same
range, although we state it only for $x>0$ real.
\end{proof}

\begin{proposition}\label{firstapr}
For $x>0$ and any natural number $N$ we have
\begin{equation}\label{e9}
\chi(x)=\sum_{n=1}^N\frac{e^{-x\sqrt{n/2}}}{\sqrt{n}}\sin\Bigl(\;\frac{\pi}{4}-x\;\sqrt\frac{n}{2}\;\Bigr)+R_N,
\end{equation}
with
\[|R_N|\le \Bigl(\frac{1}{\sqrt{N+1}}+ \frac{2\sqrt{2}}{x}\Bigr)
e^{-x\sqrt{(N+1)/2}}.\]
\end{proposition}

\begin{proof}
By \eqref{e6} we have
\[|R_n|=\Bigl|\sum_{n=N+1}^\infty\frac{e^{-x\sqrt{n/2}}}{\sqrt{n}}
\sin\Bigl(\;\frac{\pi}{4}-x\;\sqrt\frac{n}{2}\;\Bigr)\Bigr|\le \sum_{n=N+1}^\infty\frac{e^{-x\sqrt{n/2}}}{\sqrt{n}}.\]
Since the function $t\mapsto \frac{e^{-x\sqrt{t/2}}}{\sqrt{t}}$ is decreasing, we have
\[|R_N|\le \frac{e^{-x\sqrt{(N+1)/2}}}{\sqrt{N+1}}+\int_{N+1}^\infty \frac{e^{-x\sqrt{t/2}}}{\sqrt{t}}\,dt=\frac{e^{-x\sqrt{(N+1)/2}}}{\sqrt{N+1}}+
\frac{2\sqrt{2}}{x}e^{-x\sqrt{(N+1)/2}}.\]
\end{proof}

\section{An alternative integral representation.}\label{S:3}

\begin{proposition}
For any $z\in\C$ we have
\begin{equation}\label{e2}
\chi(z)=\frac{\zeta(1/2)}{\sqrt{2}}+\frac{z}{2}+2\int_0^\infty\frac{\cos(xz)-1}{e^{2\pi x^2}-1}\,dx.
\end{equation}
\end{proposition}

\begin{proof}
Let $0<\varepsilon<1$ be a given real number. 
By Cauchy's Theorem, the integration path in the definition \eqref{e1} of $\chi(z)$
can be transformed into the path consisting of the portion $(-\infty,-\varepsilon)$ 
of the negative real  axis, followed by the semicircle of center $0$ and radius 
$\varepsilon$ above the real axis, and the portion $(\varepsilon,\infty)$ of the 
real axis.
This gives us the equality
\[\chi(z)=\int_\varepsilon^\infty \frac{e^{-zxi}}{e^{2\pi x^2}-1}\,dx-
\int_{C_{\varepsilon}}\frac{e^{zwi}}{e^{2\pi w^2}-1}\,dw+
\int_\varepsilon^\infty \frac{e^{zxi}}{e^{2\pi x^2}-1}\,dx.\]
For $|w|<1$ we have
\[\frac{e^{zwi}}{e^{2\pi w^2}-1}=\frac{1}{2\pi w^2}+\frac{iz}{2\pi w}+g(w),\]
where $g(w)$ is holomorphic. Therefore,
\[\int_{C_{\varepsilon}}\frac{e^{zwi}}{e^{2\pi w^2}-1}\,dw=\frac{1}{\pi \varepsilon}
-\frac{z}{2}+G(\varepsilon),\]
with $\lim_{\varepsilon\to 0}G(\varepsilon)=0$. 
So, we obtain
\[\chi(z)=2\int_\varepsilon^\infty \frac{\cos (xz)}{e^{2\pi x^2}-1}\,dx-
\frac{1}{\pi \varepsilon}+\frac{z}{2}-G(\varepsilon).\]
We can write this in the form
\[\chi(z)=2\int_\varepsilon^\infty \frac{\cos (xz)-1}{e^{2\pi x^2}-1}\,dx+
2\int_\varepsilon^\infty \Bigl(\frac{1}{e^{2\pi x^2}-1}-\frac{1}{2\pi x^2}\Bigr)\,dx
+\frac{z}{2}-G(\varepsilon).\]
Taking the limits when $\varepsilon\to0^+$ we obtain
\[\chi(z)=2\int_0^\infty \frac{\cos (xz)-1}{e^{2\pi x^2}-1}\,dx+
2\int_0^\infty \Bigl(\frac{1}{e^{2\pi x^2}-1}-\frac{1}{2\pi x^2}\Bigr)\,dx
+\frac{z}{2}.\]
The second integral is transformed by a change of variables, and 
the resulting integral is found in Titchmarsh \cite{T}*{eq.~(2.7.1)}
\[
2\int_0^\infty \Bigl(\frac{1}{e^{2\pi x^2}-1}-\frac{1}{2\pi x^2}\Bigr)\,dx=
\frac{1}{\sqrt{2\pi}}\int_0^\infty \Bigl(\frac{1}{e^{ y}-1}-\frac{1}{y}\Bigr)\,
\frac{dy}{\sqrt{y}}=\frac{\zeta(1/2)\Gamma(1/2)}{\sqrt{2\pi}}.
\]
\end{proof}

\begin{proposition}
The function $\chi(z)$ has the following power series expansion
\begin{equation}\label{e3}
\chi(z)=\frac{z}{2}+\frac{1}{\sqrt{2}}\sum_{n=0}^\infty \frac{\zeta(n+\frac12)}{n!}\Bigl(-\frac{z^2}{8\pi}
\Bigr)^n.
\end{equation}
\end{proposition}

\begin{proof}
In the integral in equation \eqref{e2} expand the cosine in power series. It is 
easy to justify integration term by term. We obtain
\[
\chi(z)=\frac{\zeta(1/2)}{\sqrt{2}}+\frac{z}{2}+\sum_{n=1}^\infty (-1)^n \frac{z^{2n}}
{(2n)!}\cdot2\int_0^\infty \frac{x^{2n}}{e^{2\pi x^2}-1}\,dx.
\]
The integral can be computed by transforming it by a change of variables to 
the integral in Titchmarsh \cite{T}*{eq.~(2.4.1)}, and then substituting 
the known value of $\Gamma(n+\frac12)$
\begin{align*}
2\int_0^\infty \frac{x^{2n}}{e^{2\pi x^2}-1}\,dx&=\frac{2}{(2\pi)^{n+1/2}}
\int_0^\infty \frac{y^n}{e^y-1}\frac{dy}{2\sqrt{y}}=
\frac{1}{(2\pi)^{n+1/2}}
\int_0^\infty \frac{y^{n+\frac12}}{e^y-1}\frac{dy}{y}\\
&=\frac{\Gamma(n+\frac12)\zeta(n+\frac12)}{(2\pi)^{n+1/2}}=
\frac{\zeta(n+\frac12)}{(2\pi)^{n+1/2}}\frac{(2n)!\sqrt{\pi}}{n! 2^{2n}}.
\end{align*}
Therefore,
\begin{align*}
\chi(z)&=\frac{\zeta(1/2)}{\sqrt{2}}+\frac{z}{2}+\sum_{n=1}^\infty (-1)^n \frac{z^{2n}}
{(2n)!}\cdot\frac{\zeta(n+\frac12)}{(2\pi)^{n+1/2}}\frac{(2n)!\sqrt{\pi}}{n! 2^{2n}}
\\ &=\frac{\zeta(1/2)}{\sqrt{2}}+\frac{z}{2}+\frac{1}{\sqrt{2}}\sum_{n=1}^\infty
\frac{\zeta(n+\frac12)}{n!}\Bigl(-\frac{z^2}{8\pi}\Bigr)^n\\&=\frac{z}{2}+
\frac{1}{\sqrt{2}}\sum_{n=0}^\infty
\frac{\zeta(n+\frac12)}{n!}\Bigl(-\frac{z^2}{8\pi}\Bigr)^n.\qedhere
\end{align*}
\end{proof}

\begin{remark}
Riemann `derives' the power series from \eqref{e5} by means of an `Eulerian' computation,
that I cannot resist duplicating:
\begin{align*}
\chi(z)&=\frac12\sum_{n=1}^\infty \Bigl(\frac{e^{z\sqrt{n}e^{3\pi i/4}}}{\sqrt{n}}
e^{\pi i/4}
+\frac{e^{z\sqrt{n}e^{\pi i/4}}}{\sqrt{n}}e^{-3\pi i/4}\Bigr)\\
&=\frac12\sum_{n=1}^\infty \sum_{m=0}^\infty \frac{z^m n^{\frac{m-1}{2}}}{m!}\Bigl(
e^{\frac{3\pi i}{4}m+\frac{\pi i}{4}}+ e^{-\frac{3\pi i}{4}m-\frac{\pi i}{4}}\Bigr)\\&=
\sum_{n=1}^\infty \sum_{m=0}^\infty \frac{z^m n^{\frac{m-1}{2}}}{m!}\cos(\tfrac{3m+1}{4}\pi)
\\&= \sum_{m=0}^\infty \frac{z^m\zeta(\frac{1-m}{2})}{m!}\cos(\tfrac{3m+1}{4}\pi).
\end{align*}
Naturally, this is not correct, but similar calculations lead to good results on many occasions. We should not underestimate Riemann, it is quite possible that he knew a general method to justify this type of reasoning. 
The  power series obtained by Riemann is equivalent to \eqref{e3}.
\begin{equation}\label{altern}
\chi(z)=\frac{z}{2}+\frac{1}{\sqrt{2}}\sum_{n=0}^\infty \frac{\zeta(n+\frac12)}{n!}\Bigl(-\frac{z^2}{8\pi}
\Bigr)^n=\sum_{m=0}^\infty \frac{z^m\zeta(\frac{1-m}{2})}{m!}\cos(\tfrac{3m+1}{4}\pi),
\end{equation}
as can be checked easily by means of the functional equation.

Notice that we also obtain the curious equality
\begin{equation}
\chi^{(m)}(0)=\zeta\left(\tfrac{1-m}{2}\right)\cos\left(\tfrac{3m+1}{4}\pi\right).
\end{equation}
\end{remark}
\begin{corollary}
The function $\chi(z)-\frac{z}{2}$ is even so that the function $\chi(z)$ satisfies 
the functional equation
\begin{equation}\label{e4}
\chi(-z)=\chi(z)-z.
\end{equation}
\end{corollary}

\section{Connection to the zeta function.}

\begin{proposition}
We have
\begin{equation}\label{e8}
\int_0^\infty \chi(x)x^{2s-1}\,dx=\Gamma(2s)\zeta(s+\tfrac12)\cos\tfrac{\pi(2s+1)}{4},
\qquad \Re(s)>0.
\end{equation}
\end{proposition}

\begin{proof}
By \eqref{e9} with $N=1$ we have $|\chi(x)|\le (1+\frac{4}{x})e^{-x/\sqrt{2}}$ for $x>1$ and there is a constant $M$ such that $|\chi(x)|\le M$ for $0<x<1$. This implies that the left member of \eqref{e8} defines an analytic function for $\sigma>0$.  Therefore, we only have to show the equality for $\sigma=\Re(s)>1/2$.  By \eqref{e5}
\[\int_0^\infty \chi(x)x^{2s}\,\frac{dx}{x}=\int_0^\infty \frac12\sum_{n=1}^\infty \Bigl(\frac{e^{-x\sqrt{n}e^{-\pi i/4}}}{\sqrt{n}}
e^{\pi i/4}
+\frac{e^{-x\sqrt{n}e^{\pi i/4}}}{\sqrt{n}}e^{-\pi i/4}\Bigr)x^{2s}\,\frac{dx}{x}.\]
We may interchange sum and integral (for example, taking in absolute value all 
the terms and integrating we obtain a convergent series), therefore, we obtain
\begin{multline*}
\int_0^\infty \chi(x)x^{2s-1}\,dx\\=\frac12\sum_{n=1}^\infty\Bigl(\frac{e^{\pi i/4}}{\sqrt{n}}
\int_0^\infty e^{-x\sqrt{n}e^{-\pi i/4}}x^{2s-1}\,dx+
\frac{e^{-\pi i/4}}{\sqrt{n}}
\int_0^\infty e^{-x\sqrt{n}e^{\pi i/4}}x^{2s-1}\,dx\Bigr).
\end{multline*}
These integrals are elementary, and we obtain for $\sigma>1/2$
\begin{align*}
\int_0^\infty \chi(x)x^{2s-1}\,dx&=\frac12\sum_{n=1}^\infty\Bigl(\frac{e^{\pi i/4}}{\sqrt{n}}\frac{\Gamma(2s)}{n^s}e^{\frac{\pi i s}{2}}
+\frac{e^{-\pi i/4}}{\sqrt{n}}\frac{\Gamma(2s)}{n^s}e^{-\frac{\pi i s}{2}}
\Bigr)\\
&=\Gamma(2s)\zeta(s+\tfrac12 )\cos\Bigl(\frac{\pi}{4}+\frac{\pi s}{2}\Bigr).\qedhere
\end{align*}
\end{proof}

\section{Application of Euler-Maclaurin expansion.}

Equation \eqref{e5} is well suited to apply the Euler-Maclaurin expansion, which we quote from \cite{E}*{p.~104}:

\begin{proposition}\label{eulmac}
Suppose that $M$ is a positive integer and that $f$ has continuous derivatives through the $2M$th  order in the interval $[a,b]$ where $a$ and $b$ are integers  with $a<b$. Then,
\begin{multline}\label{e10}
\sum_{n=a}^b f(n)=\int_a^b f(x)\,dx+\frac{f(a)+f(b)}{2}+\sum_{m=1}^M
\frac{B_{2m}}{(2m)!}\bigl(f^{(2m-1)}(b)-f^{(2m-1)}(a)\bigr)\\
-\frac{1}{(2M)!}\int_a^b
f^{(2M)}(x)\widetilde{B}_{2M}(x)\,dx.
\end{multline}
\end{proposition}

Here, $B_m$ are the Bernoulli numbers
\[B_0=1,\quad B_1=-\frac12,\quad B_2=\frac16,\quad
B_4=-\frac{1}{30},\quad \dots,\quad   B_{2n+1}=0, \quad n\ge1\]
We also notice the inequality
\[\max_x|\widetilde{B}_{2n}(x)|=|B_{2n}|.\]

\begin{proposition}\label{p5:3}
For $x>0$ real and any natural number $N$ we have
\begin{multline}\label{ER}
\chi(x)=\sum_{n=1}^N\frac{e^{-x\sqrt{n/2}}}{\sqrt{n}}\sin\bigl(\tfrac{\pi}{4}-x\;\sqrt{n/2}\;\bigr)-\frac12\frac{e^{-x\sqrt{N/2}}}{\sqrt{N}}\sin\bigl(\tfrac{\pi}{4}-x\;\sqrt{N/2}\;\bigr)\\
-\frac{2}{x}e^{-x\sqrt{N/2}}\sin\bigl(x\;\sqrt{N/2}\;\bigr)
+\frac{x}{24N}e^{-x\sqrt{N/2}}\cos(x\sqrt{N/2})\\
+\frac{1}{24 N^{3/3}}e^{-x\sqrt{N/2}}
\cos\bigl(x\;\sqrt{N/2}+\tfrac{\pi}{4}\bigr)+R(x,N),
\end{multline}
with
\[|R(x,N)|\le \frac{e^{-x\sqrt{N/2}}}{24\sqrt{N}}(x+N^{-1/2})^2.\]
\end{proposition}

\begin{proof}
In Proposition \ref{eulmac} take  
\[f(t)= \frac12\Bigl(\frac{e^{-x\sqrt{t}e^{-\pi i/4}}}{\sqrt{t}}
e^{\pi i/4}
+\frac{e^{-x\sqrt{t}e^{\pi i/4}}}{\sqrt{t}}e^{-\pi i/4}\Bigr)=
\Re\Bigl\{\frac{e^{-x\sqrt{t}e^{-\pi i/4}}}{\sqrt{t}}e^{\pi i/4}\Bigr\}.\]
For $x>0$ we have
\[\sum_{n=N+1}^\infty f(n)=\chi(x)-\sum_{n=1}^N f(n).\]
By Euler-Maclaurin (taking $M=1$, $a=N$ and $b=L$)
\[\sum_{n=N}^L f(n)=\int_N^L f(t)\,dt+\frac{f(N)+f(L)}{2}+\frac{B_2}{2}(f'(L)-f'(N))-R(x,N,L),\]
where
\[R(x,N,L)=\frac12\int_N^L f''(t)\widetilde{B}_2(t)\,dt.\]

Taking the limit for $L\to\infty$ yields 
\[\chi(x)-\sum_{n=1}^N f(n)=\int_N^\infty f(t)\,dt-\frac{f(N)}{2}-\frac{f'(N)}{12}-\frac12\int_N^\infty f''(t)\widetilde{B}_2(t)\,dt.\]
Since
\begin{displaymath}
f''(t)=\Re\Bigl\{\frac{e^{x\sqrt{\frac{t}{2}}(-1+i)}}{4t^{3/2}}
\Bigl(\frac{3 e^{\frac{\pi}{4}i}}{t}+\frac{3x}{\sqrt{t}}+x^2e^{	\frac{7\pi}{4}i}\Bigr)
\Bigr\},
\end{displaymath}
we have
\begin{align*}
|R(x,N)|&\le \frac{1}{12}\int_N^\infty e^{-x\sqrt{t/2}}\Bigl(\frac34 t^{-5/2}+
\frac{3x}{4}t^{-2}+\frac{x^2}{4}t^{-3/2}\Bigr)\,dt\\
&\le \frac{1}{48} e^{-x\sqrt{N/2}}\bigl(2 N^{-3/2}+3 x N^{-1}+2 x^2 N^{-1/2}\bigr)\\
&\le\frac{1}{24}e^{-x\sqrt{N/2}}(N^{-3/4}+x N^{-1/4})^2.
\end{align*}
The integral of $f(t)$ can be computed and we obtain
\begin{multline*}
\chi(x)=\frac12\sum_{n=1}^N\Bigl(\frac{e^{-x\sqrt{n}e^{-\pi i/4}}}{\sqrt{n}}
e^{\pi i/4}
+\frac{e^{-x\sqrt{n}e^{\pi i/4}}}{\sqrt{n}}e^{-\pi i/4}\Bigr)
\\+
\frac{e^{-x\sqrt{N}e^{-\pi i/4}}}{x}i-
\frac{e^{-x\sqrt{N}e^{\pi i/4}}}{x}i
-\frac14 \Bigl(\frac{e^{-x\sqrt{N}e^{-\pi i/4}}}{\sqrt{N}}
e^{\pi i/4}
+\frac{e^{-x\sqrt{N}e^{\pi i/4}}}{\sqrt{N}}e^{-\pi i/4}\Bigr)\\
+\frac{1}{12} \Re\Bigl\{\frac{e^{-x\sqrt{N}e^{-\pi i/4}}}{2N}
\Bigl(\frac{e^{\frac{\pi}{4}i}}{\sqrt{N}}+x\Bigr)\Bigr\}
-R(x,N).\qedhere
\end{multline*}
\end{proof}
\begin{remark}
Riemann substitutes in \eqref{e8}  the asymptotic expansion for $\chi(s)$ obtained 
in Proposition \ref{p5:3}, integrates term by term  to get an asymptotic
expansion for $\zeta(s+\frac12)$. But he only gets the Euler-Maclaurin expansion 
of $\zeta(s+\frac12)$.
\end{remark}

Riemann uses the Euler Maclaurin expansion to compute $\chi(0)=\zeta(1/2)/\sqrt{2}$.
Which we have computed in Section   \ref{S:3}  by other means.

Taking the limits when $x\to0^+$ in \eqref{ER} he obtains:
\[\chi(0)=\sum_{n=1}^N\frac{1}{\sqrt{n}}\frac{1}{\sqrt{2}}-\frac{1}{2\sqrt{2N}}
-\sqrt{2N}
+\frac{1}{24\sqrt{2} N^{3/2}}+R(0,N).\]
Therefore,
\[\chi(0)=\lim_{N\to\infty}\Bigl(\frac{1}{\sqrt{2}}\sum_{n=1}^N\frac{1}{\sqrt{n}}-
\sqrt{2N}\Bigr).\]
Riemann recognizes the limit as $\zeta(1/2)/\sqrt{2}$.

\section{A last representation of \texorpdfstring{$\chi$}{X}.}

\begin{proposition}
For any $z\in\C$ we have
\begin{equation}\label{E:final}
\chi(z)=\frac{\zeta(\frac12)}{\sqrt{2}}+\frac{z}{2}+\sum_{n=1}^\infty\frac{e^{-\frac{z^2}{8\pi n}}-1}{\sqrt{2n}}.
\end{equation}
\end{proposition}

\begin{proof}
This can be obtained directly from the power series \eqref{e3}. Riemann does not 
write \eqref{e3} only his equivalent form in \eqref{altern}, therefore, I try 
to follow what appear to be Riemann path to this series expansion. 

Differentiating \eqref{e1} we obtain

\begin{equation}
\chi''(z)=-\int_L\frac{w^2e^{izw}}{e^{2\pi w^2}-1}\,dw.
\end{equation}
The integrand has no singularity at $w=0$, therefore we may change the line
of integration to the real line,
\[\chi''(z)=-\int_{-\infty}^\infty\frac{x^2e^{izx}}{e^{2\pi x^2}-1}\,dx.\]
From which we obtain easily
\[\chi''(z)=-\sum_{n=1}^\infty 
\int_{-\infty}^\infty x^2e^{izx}e^{-2\pi n x^2}\,dx= 
-\sum_{n=1}^\infty \frac{e^{-\frac{z^2}{8\pi n}}(4\pi n - z^2)}{16\pi^2 \sqrt{2} n^{5/2}}
.\]

\begin{figure}[H]
\begin{center}
\includegraphics[width=0.95\linewidth]{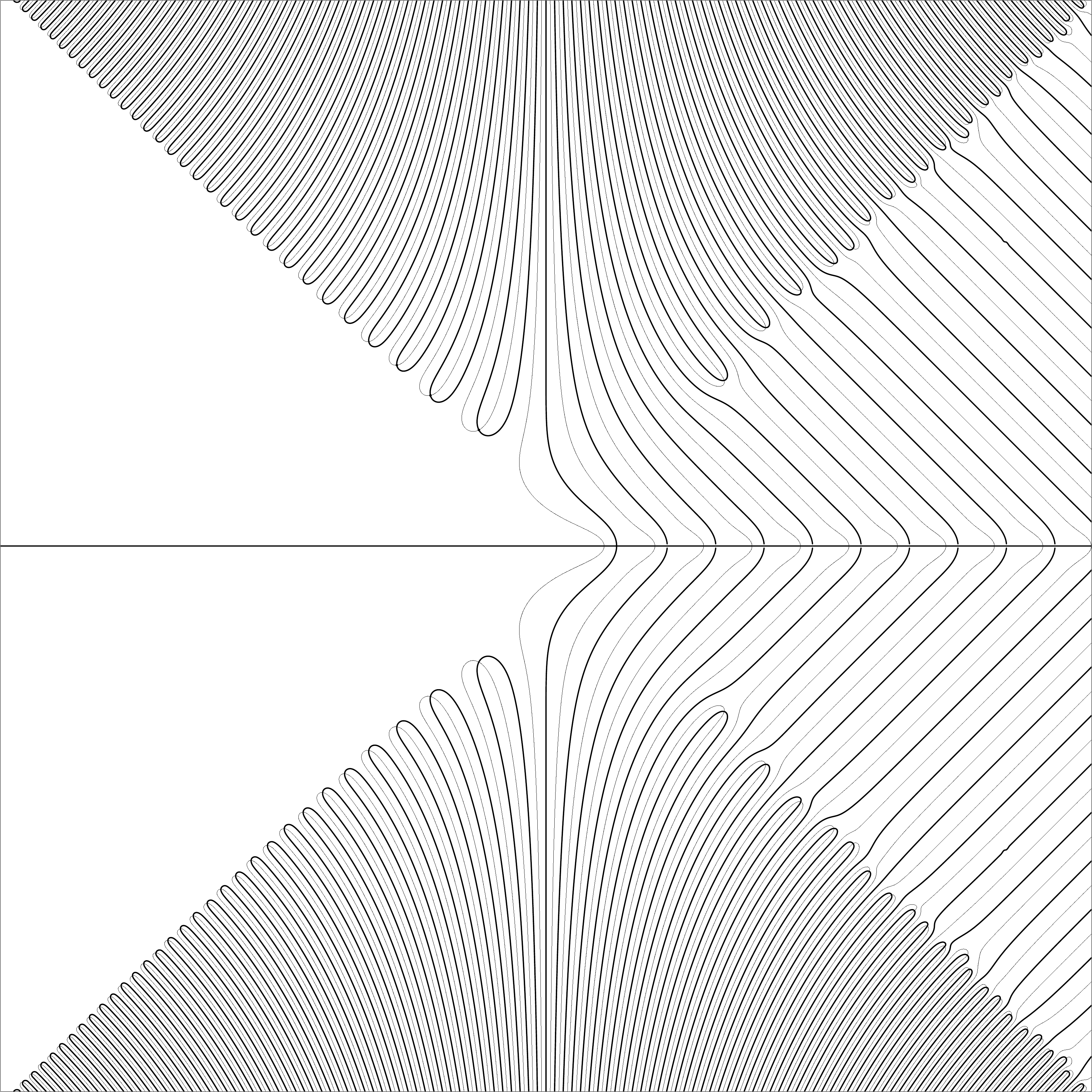}
\caption{X-ray of $\chi(s)$ on $(-50,50)^2$}
\end{center}
\end{figure}

We have 
\[\frac{\partial^2}{\partial z^2}
\frac{e^{-\frac{z^2}{8\pi n}}-1}{\sqrt{2n}}=-\frac{e^{-\frac{z^2}{8\pi n}}(4\pi n - z^2)}{16\pi^2 \sqrt{2} n^{5/2}}\]
Since the series in \eqref{E:final} converges uniformly in compact sets, we derive  
from this that 
\[\chi(z)=az+b+\sum_{n=1}^\infty\frac{e^{-\frac{z^2}{8\pi n}}-1}{\sqrt{2n}},\]
for some adequate constants $a$ and $b$. Taking $z=0$, we obtain $b=\chi(0)$ and
it is also clear that $\chi'(0)=a$.  This proves \eqref{E:final}.
\end{proof}

\section{x-ray and zeros.}

We may use the power series to compute $\chi(z)$ for $|z|\le r$. Since there
is a large cancellation between the terms of the power series, we have to 
use large precision in computing the terms. Thus, we have been able 
to compute $\chi(z)$ in the square $(-50,50)^2$. This allows us to
compute the x-ray of this function in this square.

Taking $N=1$ in Proposition \ref{firstapr} we have for $x>0$
\[\chi(x)\sim e^{-x/\sqrt{2}}\sin\Bigl(\frac{\pi}{4}-\frac{x}{\sqrt{2}}\Bigr)+
O(e^{-x}).\]
Therefore, the real zeros are approximately equal to $a_n\sim \pi\sqrt{2}(n+\frac14)$.
In practice, this is a very good approximation; for example, for $n=10$ the difference
$a_{10}$ and its approximation is $1.12\times10^{-6}$. 

The zeros in the direction $e^{3\pi i/4}$ are approximately
situated at $4\pi\sqrt{n}e^{3\pi i/4}$. But this is only a guess founded on the numbers,
and it is not so good as in the real case. 

The zeros in the direction $e^{\pi i/4}$ are more difficult to predict because
some of them \emph{are missing} because they are \emph{migrated} to the real line.

\end{document}